\documentclass[12pt]{amsart}
\usepackage{amsrefs} 
\usepackage{calc,amssymb,amsmath,amsthm,amsfonts,ifthen,graphics,enumerate}
\usepackage[all]{xy}
\usepackage[ansinew]{inputenc} 
\usepackage{accents}
\usepackage{epsfig,amssymb,amsmath,amsthm,enumerate}
\usepackage{float}
\newfloat{xyfloat}{htb}{xyfl}
\usepackage{wrapfig}

\newcommand{\topplace}{}

\newcommand{\ignore}[1]{\relax}

\renewcommand{\date}[1]{\pagebreak[3]\marginpar{#1}\nopagebreak[4]}

\newcommand{\defem}[1]{%
  \ifmmode{\ #1\ }\else{{\mdseries\itshape\sffamily #1}}\fi\index{#1}}

\newcommand{\ident}{\ensuremath{\approx}}

\newcommand{\Assoc}{\text{Assoc}}

\newcommand{\mapsfrom}{\kern1.5pt\leftarrow\joinrel\mapsfromchar}

\newcommand{\mapsfromchar}{%
\begin{picture}(1,5)(2,-0.75)
\put(0,0){\line(0,1){4.5}}
\end{picture}}

\renewcommand{\iff}{\text{ iff }}

\newcounter{tabitem}

\newcounter{delta}
\newcounter{doubledelta}
\newcounter{yA0}
\newcounter{yA1}
\newcounter{yA2}
\newcounter{yA3}
\newcounter{yA4}
\newcounter{yA5}
\newcounter{yB0}
\newcounter{yB1}
\newcounter{yB2}
\newcounter{yB3}
\newcounter{yB4}
\newcounter{yB5}

\newcounter{xA0}
\newcounter{xA1}
\newcounter{xA2}
\newcounter{xA3}
\newcounter{xA4}
\newcounter{xA5}
\newcounter{xB0}
\newcounter{xB1}
\newcounter{xB2}
\newcounter{xB3}
\newcounter{xB4}
\newcounter{xB5}

\newcounter{locx1}
\newcounter{locx2}
\newcounter{locy1}
\newcounter{locy2}

\newcommand{\blob}{\circle*{5}}
\newcommand{\vertex}[1]{\put #1{\blob}}
\newcommand{\edge}[2]{\qbezier#1 #2 #2}

\newcounter{poly}
\newcommand{\xpoly}[4]{\setcounter{poly}{#1}
  \vertex{#3}
  \edge{#3}{#4}
  \ifthenelse{\value{poly}>3}
	     {\addtocounter{poly}{-1}\xpoly{\value{poly}}{#2}{#4}}%
	     {\vertex{#4}
	      \edge{#4}{#2}}
}

\newcommand{\circuit}[3]{\setcounter{poly}{#1}
  \vertex{#2}
  \vertex{#3}
  \edge{#2}{#3}
  \ifthenelse{\value{poly}>2}{
    \xpoly{#1}{#2}{#3}}{}
}

\newcommand{\poly}[3]{\setcounter{poly}{#1}
  \vertex{#2}
  \edge{#2}{#3}
  \ifthenelse{\value{poly}>2}
	     {\addtocounter{poly}{-1}\poly{\value{poly}}{#3}}%
	     {\vertex{#3}}
}

\newcounter{rot}

\newsavebox{\cbox}
\newlength{\boxwidth}
\newlength{\boxheight}
\newlength{\boxdepth}

\newcommand{\hcbox}[1]{\savebox{\cbox}{\mbox{#1}}%
  \settowidth{\boxwidth}{\usebox{\cbox}}%
  \setlength{\boxwidth}{0.5\boxwidth}%
  \mbox{\kern-\boxwidth\usebox{\cbox}}}

\newcommand{\hcfbox}[1]{\savebox{\cbox}{\fbox{#1}}%
  \settowidth{\boxwidth}{\usebox{\cbox}}%
  \setlength{\boxwidth}{0.5\boxwidth}%
  \mbox{\kern-\boxwidth\usebox{\cbox}}}

\newsavebox{\lmp}
\newcommand{\leftmarginpar}[1]{\savebox{\lmp}{\parbox{1in}{#1}}%
\ \kern-1in\kern-\leftmargin\raisebox{\baselineskip}{\usebox{\lmp}}}

\newcommand{\newthm}[3]{\ifthenelse{\equal{#3}{}}{%
\newtheorem{#1}{#2}}{%
\newtheorem{#1}{#2}[#3]}}

\newcommand{\thmplace}{\topplace}
\theoremstyle{plain}
\newthm{theorem}{Theorem}{\thmplace}

\newtheorem{proposition}[theorem]{Proposition}
\newtheorem{lemma}[theorem]{Lemma}
\newthm{axiom}{Axiom}{\thmplace}
\newthm{conjecture}{Conjecture}{\thmplace}
\newthm{algorithm}{Algorithm}{\thmplace}
\theoremstyle{definition}
\newthm{definition}{Definition}{\thmplace}

\newthm{example}{Example}{\thmplace}

\newthm{Note}{Note}{\thmplace}

\newtheorem*{ntheorem}{Theorem}

\newcounter{alfa}[theorem]
\renewcommand{\thealfa}{\alph{alfa})}

\newcounter{theline}

\newlength{\caretsize}
\newlength{\strhgt}
\newcommand{\str}[2][\strhgt]{\raisebox{-#1}{$#2$}}

\newcommand{\caret}[5][{}]{
\xy <0mm,0mm>;<\caretsize,0mm>:
(-#2,-#3)="l";  
(0,0)="t";   
(#2,-#3)="r";   
"l";"t"**@{-};
    "r"**@{-};
"l"*!UC=(#2,#3){#4};
"r"*!UC=(#2,#3){#5};
"t"*{\raisebox{-3pt}[0pt][0pt]{#1}};
\endxy}

\newcommand{\ldcaret}[5][{}]{
\xy <0mm,0mm>;<\caretsize,0mm>:
(-#2,-#3)="l";  
(0,0)="t";   
(#2,-#3)="r";   
"l";"t"**@{.};
    "r"**@{-};
"l"*!UC=(#2,#3){#4};
"r"*!UC=(#2,#3){#5};
"t"*{\raisebox{-3pt}[0pt][0pt]{#1}};
\endxy}

\newlength{\trdp}
\newlength{\trht}
\newlength{\trrs}
\newlength{\trlf}
\newlength{\trrt}
\newcommand{\ltree}[4][0pt]{%
\setlength{\trrs}{25pt}\addtolength{\trrs}{#1}
\setlength{\trht}{15pt}\addtolength{\trht}{#1}
\setlength{\trdp}{22pt}\addtolength{\trdp}{#1}
\setlength{\trlf}{0pt}\addtolength{\trlf}{#1}
\setlength{\trrt}{-0pt}\addtolength{\trrt}{#1}
\raisebox{\trrs}[\trht][\trdp]{\hspace{\trlf}
{\caret{5}{10}{\caret{5}{10}{\str{#2}}{\str{#3}}}{\str{#4}}}\hspace{\trrt}}
}

\newcommand{\rtree}[4][0pt]{%
\setlength{\trrs}{25pt}\addtolength{\trrs}{#1}
\setlength{\trht}{15pt}\addtolength{\trht}{#1}
\setlength{\trdp}{22pt}\addtolength{\trdp}{#1}
\setlength{\trlf}{-3pt}\addtolength{\trlf}{#1}
\setlength{\trrt}{3pt}\addtolength{\trrt}{#1}
\raisebox{\trrs}[\trht][\trdp]{\hspace{\trlf}
{\caret{5}{10}{\str{#2}}{\caret{5}{10}{\str{#3}}{\str{#4}}}}\hspace{\trrt}}
}

\newcommand{\carethwidth}{5}
\newcommand{\caretheight}{10}
\newcommand{\onecaret}[3][15pt]{\raisebox{#1}{\caret{\carethwidth}{\caretheight}{\str{#2}}{\str{#3}}}}
\newcommand{\onecaretw}[3][15pt]{\raisebox{#1}{\caret{10}{\caretheight}{\str{#2}}{\str{#3}}}}

\newcommand{\treecolor}{\relax}
\newcommand{\formulacolor}{\relax}

\newcommand{\ovln}[1]{\overline{#1}}
\newcommand{\assoc}{\textup{Assoc}}

\newcommand{\beginproof}[1][\relax]%
{\ifthenelse{\equal{\relax}{#1}}
{\noindent\textit{\proofname.\ }}{\noindent\textit{#1.\ }}}
\renewcommand{\endproof}{\qed}

\begin{document}

\title{Associativity of the Commutator Operation in Groups}
\author[F. Guzmán]{Fernando Guzmán}
\address{{\flushleft Fernando Guzmán}
\newline\indent
Department of Mathematical Sciences
\newline\indent
Binghamton University
\newline\indent
Binghamton N.Y., 13902-6000}
\email{fer@math.binghamton.edu}
\urladdr{http://math.binghamton.edu/fer}


\begin{abstract}
The study of associativity of the commutator operation in
groups goes back to some work of Levi in 1942.  In the 1960's Richard
J. Thompson created a group F whose elements are representatives of
the generalized associative law for an arbitrary binary operation.  In
2006, Geoghegan and Guzmán proved that a group G is solvable if and only if the
commutator operation in G eventually satisfies ALL instances of the
associative law, and also showed that many non-solvable groups do not
satisfy any instance of the generalized associative law.  We will
address the question: Is there a non-solvable group which satisfies
SOME instance of the generalized associative law?  For finite groups,
we prove that the answer is no.
\end{abstract}

\maketitle

\section{Introduction}

In 1942 F.~W.~Levi \cite{levi} proved that the commutator operation in
a group is associative if and only if the group is nilpotent of class
$\leq 2$. That, in a sense, settled the question of associativity of
the commutator operation in groups.  In the 1960's Richard Thompson,
studying the logical connections between the three-variable
associative law and the generalized associative law, came up with the
group $F$ whose elements can be thought of as (equivalence classes of)
instances of the generalized associative law.  One way to encode a
parenthezation of a product of $n$ factors, is using a binary tree
with $n$ leaves; that way, each instance of the generalized
associative law, is encoded with a pair of binary trees having the
same number of leaves.  The equivalence relation that gives rise to
the group $F$ declares two pairs of trees to be equivalent if one pair
is obtained from the other by doing a ``{\sl common expansion}'' on
both trees. The precise details can be found in the paper by Cannon,
Floyd and Parry \cite{cfp}.  We illustrate the idea with the following
example. If $t_1,t_2,t_3$ are any binary trees then the following two
pairs are equivalent.

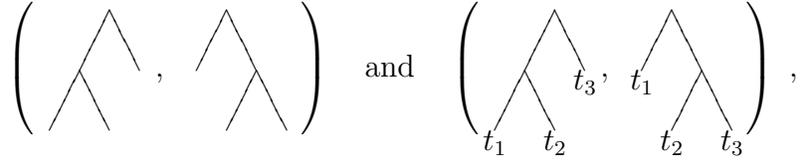
\[ \left(\ \ltree{}{}{},\rtree{}{}{}\ \ \right) \text{\quad and\quad }
   \left(\ \ltree{t_1}{t_2}{t_3},\rtree{t_1}{t_2}{t_3}\ \ \right)\ , \]

\vspace{10pt}
\noindent
as the second pair is obtained from the first pair by ``{\sl hanging}''
$t_1$ from the first leaf, $t_2$ from the second leaf and $t_3$ from
the third leaf.  Note that the left pair encodes the 3-variable
associative law.  For each equivalence class, there is a ``{\sl
  reduced pair}'' of trees that represent it, that is, one pair from
which all the others in that equivalence class can be obtained by some
common expansion.  A reduced pair is characterized by the fact that
the two trees do not have any ``{\sl matching}'' free caret.  The left
pair above is reduced.

In 2006 Geoghegan and Guzmán \cite{geog-guzm} used subgroups of
Thompson's group $F$ to ``{\sl measure}'' the associativity of a
binary operation.   That is, a non-associative binary operation may
still satisfy some instances of the generalized associative law.
In particular, a group which is not nilpotent of class $\leq 2$, may
satisfy some instances of the generalized associative law, and this
reopens the issue of the associativity of the commutator operation in
groups, in this generalized sense.  
Because the elements of $F$ are equivalence classes, we say that a
binary operation ``{\sl eventually satisfies}'' an instance of the
generalized associative law, if it satisfies some expansion of it,
that is, an instance in the same equivalence class.  For any magma
$S$, the set $\assoc(S)$ of instances of the generalized associative
law that $S$ satisfies, is a subgroup of $F$, and it is proved in
\cite{geog-guzm} that when we take $S$ to be a group $G$ with its
commutator operation, $G$ is a solvable group if and only if $\assoc(G)=F$.
Moreover, for many non-solvable groups $G$ it was shown that
$\assoc(G)=1$.  That opened the question of whether there is a
non-solvable group $G$ for which 
  \[ 1 \lneq \Assoc(G) \lneq F  \]
i.e. such that the commutator operation of $G$ satisfies some instance
of the generalized associative law.  

In this paper, we prove

\begin{ntheorem}\label{thm:intro}
  If the commutator operation a finite group $G$ satisfies some
  instance of the generalized  associative law, then $G$ is solvable.
\end{ntheorem}

\section{Notation}

Let $G$ be a group $x,y\in G$. We will use the standard notation
\[ [x,y] := x^{-1}y^{-1}xy \]
for the commutator of $x$ and $y$. 
We  will use binary trees to denote
general commutator expressions, according to the following recursive
definition.  If $s$ and $t$ are commutator expressions, then the
binary tree 
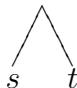
\[ \caret{5}{10}{\str[8pt]{s}}{\str[8pt]{t}} \]
denotes the commutator expression $[s,t]$.  This way, the $3$-variable
associative law for the commutator operation in a group can be written
as
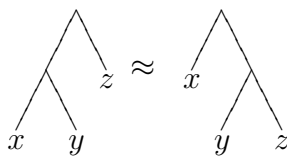
\[ \ltree{x}{y}{z} \ident \rtree{x}{y}{z} \]
\vspace{20pt}

A \defem{tree-like} expression in a group $G$ is a commutator
expression, where each leaf of the corresponding binary tree is
labelled with either a variable, or an element of the group $G$.  We
will refer to the elements of $G$ in a tree-like expression as
\defem{constants}.  In particular, if we take the trivial binary tree,
the one having a single leaf, then any element of $G$ is a constant
tree-like expression in $G$.  

  Given two tree-like expressions $s$ and
$t$ in a group $G$, and a subset $X\subseteq G$, we say that $X$
satisfies the identity $s\ident t$, and write $X\models s\ident t$, if
we get equality whenever we substitute elements of $H$ for the
variables of $s$ and $t$.  Given a tree-like expression $t$ in a group
$G$, and a subset $X\subseteq G$, we denote by $t(X)$ the set of all
values obtained by substituting elements of $X$ for the variables of
$t$.

A subset $X\subseteq G$ of a group is said to be \defem{normal} if it
is closed under conjugation by every element of $G$, and it is said to
be \defem{inverse} if it contains the inverse of each of its elements.
We denote by $B_p$ the full binary tree of height $p$.  The set
$B_p(G)$ is a normal, inverse, generating set of the $p$-th derived
group $G^{(p)}$.

\section{Basic Results}

The commutator operation satisfies a number of identities in all
groups, some of which we need to refer to explicitly.

\begin{eqnarray}
\label{eqn:comm 1}
  [xy,z] &=& [x,z]^y\cdot[y,z]  \\
\label{eqn:comm 2}
  [x,yz] &=& [x,z]\cdot[x,y]^z  \\
\label{eqn:comm 3}
  [y,x] &=& [x,y]^{-1} = [x^y,y^{-1}] = [x^{-1},y^x] 
\end{eqnarray}

If $G$ is a group and $s,t,s',t'$ are tree-like expressions such that 
 $(s',t')$ is an expansion of $(s,t)$, then  $G\models s\ident
t$ clearly implies that $G\models s'\ident t'$.  The converse does not
hold, but if the expansion from $(s,t)$ to $(s',t')$ involves hanging
trees of height at most $p$ then if $G\models s'\ident t'$ it follows
that $B_p(G)\models s\ident t$.  We can take this one step further
when one side of the identity is trivial.

\begin{lemma}\label{lemma:passing to G^p} 
  Let $G$ be a group,  $t$ a tree-like expression and $t'$ the
  result of hanging from each non-constant leaf of $t$ a tree of 
  height at most $p$. \\  
  If {$G\models t'\ident 1$}  then $B_p(G)\models t\ident
  1$ and {$G^{(p)}\models t\ident 1$}.
\end{lemma}

\begin{proof}
  Since every tree of height at most $p$ can be expanded to a full
  binary tree $B_p$ of height $p$, 
  from $G\models t'\ident 1$ we immediately get $B_p(G)\models t\ident 1$.
  Now, the  commutator  identities~(\ref{eqn:comm 1}) and~(\ref{eqn:comm 2})
  and  the fact that $B_p(G)$ is a normal, inverse,
  generating set of $G^{(p)}$, yield $G^{(p)}\models t\ident 1$.
\end{proof}

The next lemma follows by a straight forward calculation. When
evaluating $t$ and $s$ in $G/Z(G)$ the constants are to be replaced by
their corresponding cosets.

\begin{lemma}\label{lemma:mod the center}
  Let $G$ be a group, $Z(G)$ its center, and let $s$, $t$ be tree-like
  expressions in $G$. Let $x$ be a variable, different from all
  variables that occur in $s$ and $t$. \\
  {\small  \setlength{\strhgt}{9pt}
  \formulacolor{$G/Z(G)\models s\ident t$} \iff
  \formulacolor{$G\models \onecaret{{s}}{{x}} \ident
    \onecaret{{t}}{{x}}$} 
\iff
  \formulacolor{$G\models \onecaret{{x}}{{s}} \ident \onecaret{{x}}{{t}} $}}
\end{lemma}

A \defem{free caret} in a binary tree, is a caret with two leaves.
Each non-trivial binary tree has at least one caret; if it has exactly one caret,
we call it a \defem{vine}.  The \defem{left vine} of height $n$,
denoted $l_n$ is the vine of height $n$, whose free caret holds the
two leftmost leaves of the vine.  In the following picture we see a
vine $v_5$ of height 5 and the left vine $l_5$.

\[\treecolor{  v_5: \raisebox{45pt}{
            \caret{5}{6}{}{
              \caret{5}{6}{
		\caret{5}{6}{}{
		  \caret{5}{6}{}{
		    \caret{5}{6}{}{}}}}{}}}
\hspace{80pt}
  l_5:\hspace{30pt} \raisebox{45pt}{
            \caret{5}{6}{
              \caret{5}{6}{
		\caret{5}{6}{
		  \caret{5}{6}{
		    \caret{5}{6}{}{}}{}}{}}{}}{}}
}
\]

\vspace{30pt}

Given a vine $v_n$ of height $n$, $a\in G$, $u=(x_1,x_2,\dots,x_n)$ we
denote by $v_{n.l}(a,u)$ (resp. $v_{n.r}(a,u)$) the tree like
expression in $G$ obtained by placing $a$ at the left (resp. right)
leaf of the free caret in $v_n$ and $x_1,\dots,x_n$ at the other
leaves of $v_n$ from bottom to top.  In the previous example we have

\[
  v_{5,l}(a,u): \raisebox{45pt}{
            \caret{5}{6}{\str{x_5}}{
              \caret{5}{6}{
		\caret{5}{6}{\str{x_3}}{
		  \caret{5}{6}{\str{x_2}}{
		    \caret{5}{6}{\str{a}}{\str{x_1}}}}}{\str{x_4}}}}
\hspace{80pt}
  {v_{5,r}}(a,u): \raisebox{45pt}{
            \caret{5}{6}{\str{x_5}}{
              \caret{5}{6}{
		\caret{5}{6}{\str{x_3}}{
		  \caret{5}{6}{\str{x_2}}{
		    \caret{5}{6}{\str{x_1}}{\str{a}}}}}{\str{x_4}}}}
\]

\vspace{30pt}

\begin{lemma}\label{lemma:left vine}
  Let $v_n$ be a vine of height $n$, $G$ a group, $a\in G$, and
  $u=(x_1,x_2,\dots,x_n)$.   Then
  \[ v_{n,l}(a,u)=l_{n,l}(\ovln{a},\ovln{u}) = l_{n,l}(a,\widehat{u})^{\pm 1} \]
  where $\ovln{a}\in\{a,a^{-1}\}$ and each $\ovln{x_i},\widehat{x_i}$ are conjugates of
  $x_i$.  Similarly,
  \[ v_{n,r}(a,u)=l_{n,l}(\ovln{a},\ovln{u}) = l_{n,l}(a,\widehat{u})^{\pm 1} \]
  where $\ovln{a}\in\{a,a^{-1}\}$ and $\ovln{x_i},\widehat{x_i}$ are conjugates of
  $x_i$.
\end{lemma}

\begin{proof}
  When $n=1$, it follows from the commutator identities~(\ref{eqn:comm
    3}).\\ 
  When $n>1$, it follows from those same identities and induction.
  There are four cases to consider, two for $v_{n,l}(a,u)$ and two for
  $v_{n,r}(a,u)$.  If 
  \[
  v_{n,l}(a,x_1,\dots,x_n)=[x_n,v_{n-1,l}(a,x_1,\dots,x_{n-1})]=
  \onecaretw{x_n}{v_{n-1,l}(a,u)}, \]
  by induction, 
  $v_{n-1,l}(a,x_1,\dots,x_{n-1})^{-1}=l_{n-1,l}(\ovln{a},\ovln{x_1},\dots,\ovln{x_{n-1}})$ where
  $\ovln{a}\in\{a,a^{-1}\}$ and each $\ovln{x_i}$ is a conjugate of $x_i$,
  $i=1,\dots n-1$.  Now, using~(\ref{eqn:comm 3}), we get
  \[ 
  \begin{array}{rclcl}
    v_{n,l}(a,x_1,\dots,x_n)
    &=&[x_n,l_{n-1,l}(\ovln{a},\ovln{x_1},\dots,\ovln{x_{n-1}})^{-1}] 
    &=&\mkern-35mu\onecaretw{x_n}{l_{n-1,l}(\ovln{a},\ovln{u})^{-1}} \\
  &=&[l_{n-1,l}(\ovln{a},\ovln{x_1},\dots,\ovln{x_{n-1}}),\ovln{x_n}] 
    &=&{\setlength{\strhgt}{9pt}
      \onecaretw{l_{n-1,l}(\ovln{a},\ovln{u})}{\ovln{x_n}}} \\
  &=&l_{n,l}(\ovln{a},\ovln{x_1},\dots,\ovln{x_{n-1}},\ovln{x_n})
  \end{array}
  \]
  where $\ovln{x_n}=(x_n)^{l_{n-1,l}(a',x_1',\dots,x_{n-1}')^{-1}}$. 
  The other case for $v_{l,n}$
  \[
  v_{n,l}(a,x_1,\dots,x_n)=[v_{n-1,l}(a,x_1,\dots,x_{n-1}),x_n]=
  \mkern20mu\onecaretw{v_{n-1,l}(a,u)}{x_n}, \]
  and the two cases for $v_{n,r}$ are treated similarly.
  The second equality of the statement, i.e. 
  \[ l_{n,l}(\ovln{a},\ovln{u}) = l_{n,l}(a,\widehat{u})^{\pm 1} \]
  is obtained using~(\ref{eqn:comm 3}) repeatedly to ``{\sl pull}'' the
  inverse out of the commutator when the exponent of $a$ is $\ ^{-1}$.
\end{proof}

We illustrate the first part of this lemma with the following example:

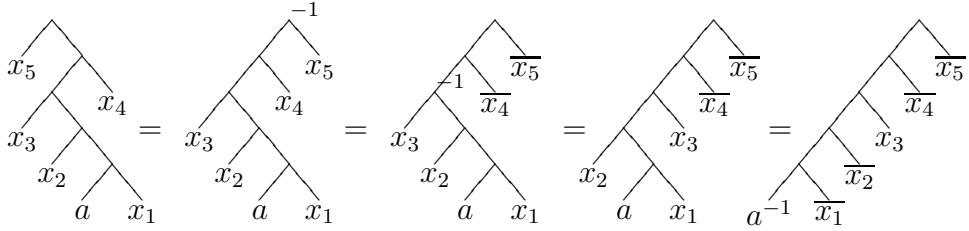
\[\hspace{-10pt}
 \raisebox{45pt}{
            \caret{5}{6}{\str{x_5}}{
              \caret{5}{6}{
		\caret{5}{6}{\str{x_3}}{
		  \caret{5}{6}{\str{x_2}}{
		    \caret{5}{6}{\str{a}}{\str{x_1}}}}}{\str{x_4}}}}
\hspace{12pt}=\hspace{20pt}
 \raisebox{45pt}{
            \caret[\mbox{\ \ $\ ^{-1}$}]{5}{6}{
              \caret{5}{6}{
		\caret{5}{6}{\str{x_3}}{
		  \caret{5}{6}{\str{x_2}}{
		    \caret{5}{6}{\str{a}}{\str{x_1}}}}}{\str{x_4}}}{\str{x_5}}}
\hspace{0pt}=\hspace{20pt}
 \raisebox{45pt}{
            \caret{5}{6}{
              \caret{5}{6}{
		\caret[\mbox{\ \ $\ ^{-1}$}]{5}{6}{\str{x_3}}{
		  \caret{5}{6}{\str{x_2}}{
		    \caret{5}{6}{\str{a}}{\str{x_1}}}}}{\str{\ovln{x_4}}}}{\str{\ovln{x_5}}}}
\hspace{5pt}=\hspace{20pt}
 \raisebox{45pt}{
            \caret{5}{6}{
              \caret{5}{6}{
		\caret{5}{6}{
		  \caret{5}{6}{\str{x_2}}{
		    \caret{5}{6}{\str{a}}{\str{x_1}}}}{\str{x_3}}}{\str{\ovln{x_4}}}}{\str{\ovln{x_5}}}}
\hspace{0pt}=\hspace{20pt}
 \raisebox{45pt}{
            \caret{5}{6}{
              \caret{5}{6}{
		\caret{5}{6}{
		  \caret{5}{6}{
		    \caret{5}{6}{\str{a^{-1}}}{\str{\ovln{x_1}}}}{\str{\ovln{x_2}}}}{\str{x_3}}}{\str{\ovln{x_4}}}}{\str{\ovln{x_5}}}}
\]

\vspace{30pt}

As an immediate consequence of Lemma~\ref{lemma:left vine} we obtain:

\begin{proposition}\label{prop:centralize}
  Let $G$ be a group, $a,b\in G$, and $v_n$ a vine of height $n$. \\
  If $b$ centralizes $l_{n,l}(a,u)$ for 
  all $u$, \\
  then it also centralizes $v_{n,l}(a,u)$ and  ${v_{n,r}}(a,u)$
  for all $u$. 
\end{proposition}

\section{The Main Result}

Before we prove the main result of this paper, we will need a counting
argument which can be expressed in terms of coloring of the leaves of
a full binary tree.

Given two leaves in a full binary tree, we'll refer to the distance to
their closest common ancestor, as the \defem{distance} between the
leaves.

\begin{proposition}\label{prop:coloring}
  Suppose we  color the leaves of a full binary tree of height $nj+1$
  so that the coloring satisfies the following condition: 
  \textrm{any two leaves at a distance }
  \[  \formulacolor{d\equiv 1 \pmod{j}}  \]
  \textrm{ must have different color.}
  Then the  number of colors has to be $\geq 2^n$.
\end{proposition}

\begin{proof}
  The proof is by induction on $n$.  When $n=1$, any leaf on the left
  subtree is at a distance $j+1$ from any leaf on the right subtree.
  Hence at least two colors are needed.  For $n>1$, consider the $2^j$
  subtrees, $t_1,\dots,t_{2^j}$ which are full binary of height
  $(n-1)j+1$.  Each one of 
  them needs at least $2^{n-1}$ colors.  Each of the leaves on the
  leftmost subtree, $t_1$, are at a distance $nj+1$ from each of the
  leaves on the rightmost subtree, $t_{2^j}$, hence all the $2^{n-1}$
  colors used for $t_1$ must be different from the $2^n$ used for
  $t_{2^j}$, and we need $2^n$ colors.
\end{proof}

Although we don't need it here, it is not hard to see that this lower
bound is tight, i.e. one can always do the coloring with $2^n$ colors.

We now get to the main result of this paper.

\begin{theorem}\label{thm:main}
  If the commutator operation a finite group $G$ satisfies some
  instance of the generalized  associative law, then $G$ is solvable.
\end{theorem}

\beginproof
  Let $s',t'$ be two binary trees with the same number of leaves, so
\begin{wrapfloat}{xyfloat}{i}[0.2\width]{1in}
\raisebox{80pt}[80pt][100pt]{
\xy <0mm,0mm>;<\caretsize,0mm>:
(-25,-20)*{t:};
(0,0)*{\bullet}; (-10,-40)*{k}
**\crv{~*=<8pt>{.} (-6,-10)&(6,-20)&(0,-30)};
(-10,-52)*{\caret{5}{10}{x_i}{
         \caret{5}{10}{
 	  \ldcaret{5}{10}{
 	    \caret{5}{10}{x_{i+1}}{}}{}}{}}}
\endxy}
\end{wrapfloat}
  that  $s'\ident t'$ is an instance of the generalized associative
  law. Let $G$ be a finite group such that $G\models s'\ident t'$.  Let
  $(s,t)$ be a reduced pair of trees, 
  which represents the equivalence class $[s',t']\in F$.   Since
  $(s',t')$ is obtained from 
  $(s,t)$ by a common expansion, there is a $p\geq 0$ such that
  $B_p(G)\models s\ident t$.  In fact, $p$ can be taken to be the
  max of the heights of the trees used to expand $(s,t)$ into
  $(s',t')$.  Label the variables of both $s$ and $t$ with variable
  names $x_1,x_2,\dots$ from left to right, and let
  {\renewcommand{\caretheight}{5} $\onecaret[10pt]{x_i}{x_{i+1}}$}\ \ be the
  leftmost free caret in either $s$ or $t$; without lost of
  generality, let's say in $s$.  Let $k$ be the 
  lowest common ancestor of $x_i$ and $x_{i+1}$ in $t$, and $r$ the
  subtree of $t$ rooted at $k$.
  Since $t$ has no free caret to the left of $x_{i+1}$, the left
  child of $r$ has no free caret, and hence it has to be the leaf
  $x_i$; $x_{i+1}$ is the first leaf of the right child of $r$. So $r$
  is an expansion of
  {\renewcommand{\caretheight}{5}
    $\mkern-5mu\onecaretw[20pt]{x_i}{l_{j,l}(x_{i+1},y)}$}
  \ \quad
  for variables $y_1,\dots,y_j$. Consider the path from the root of $t$ to
  the leaf $x_{i+1}$.  This path goes through the vertex $k$, and
  determines a vine $v_{m,l}(x_i,l_{j,l}(x_{i+1},y),z_2,\dots,z_m)$,
  so that $t$ is an expansion of this vine. Let
  $p'$ be the max of the heights of the trees hanged at 
  $y_1,\dots,y_j,z_2,\dots,z_m$ to get $t$. For any $a$ and $b$ that
  commute, putting them in the place of variables $x_i$ and
  $x_{i+1}$, evaluates $s$ to $1$, and therefore, by
  Lemma~\ref{lemma:passing to G^p}
  \[ B_{p+p'}(G)\models v_{m,l}(b,l_{j,l}(a,y),z_2,\dots,z_m) \ident 1 \]
  and
  \[ G^{(p+p')}\models v_{m,l}(b,l_{j,l}(a,y),z_2,\dots,z_m) \ident 1 \]
  Repeated application of Lemma~\ref{lemma:mod the center} yields
  \[ G^{(p+p')}/Z_{m-1}(G^{(p+p')})\models [b,l_{j,l}(a,y)]. \ident 1 \]
  Let $H=G^{(p+p')}/Z_{m-1}(G^{(p+p')})$.  Thus, we have shown the
  following fact for $H$: if $b$ commutes with $a$ then $b$
  centralizes $l_{j,l}(a,u)$ for all $u_i\in H$.  Applying this fact
  to $b$ and  $l_{j,l}(a,u)$, we conclude that $b$ centralizes
  $l_{2j,l}(a,u)$, and by induction on $q$ we get that $b$ centralizes
  $l_{qj,l}(a,u)$.  By Proposition~\ref{prop:centralize}, we get that
  for any vine $v_{qj}$, $b$ centralizes $v_{qj,l}(a,u)$ and
  $v_{qj,r}(a,u)$.  Now, applying the fact to $v_{qj,l/r}(a,u)$ and
  $b$, induction, and Proposition~\ref{prop:centralize}, we conclude
  that $v_{qj,l/r}(a,u)$ commutes with $w_{qj,l/r}(b,u')$ for any
  vines $v_{qj}$ and $w_{qj}$, and any $u,u'\in H$.  So, in a full
  binary tree of height $qj+1$. if two of the leaves, one in the left
  subtree and one in the right subtree, are labelled with
  elements that commute, the whole tree evaluates to $1$.

\ \vspace{30pt}
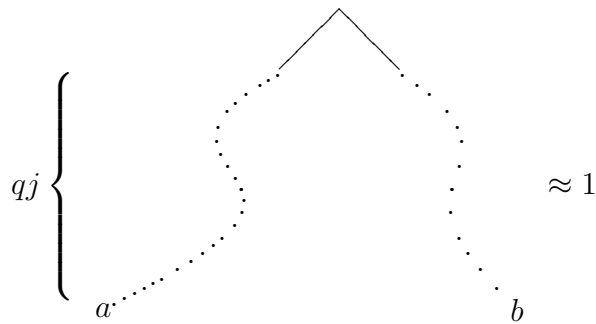
\[
qj \left\{\hspace{35pt}
\rule[-40pt]{0pt}{80pt}\quad\quad
\raisebox{70pt}[30pt][20pt]{
\caret{10}{10}{
\xy <0mm,0mm>;<\caretsize,0mm>:
(30,0)*{};
(-1,-1)*!UC=(0,0){}; (-30,-40)*{a}
**\crv{~**=<8pt>{.} (-18,-10)&(0,-20)&(-10,-30)};
\endxy}{
\xy <0mm,0mm>;<\caretsize,0mm>:
(-20,0)*{};
(1,-1)*!UC=(0,0){}; (20,-40)*{b}
**\crv{~*=<8pt>{.} (15,-10)&(6,-20)&(10,-30)};
\endxy}
}
\right. \hspace{40pt}\ident 1
\]
\vspace{10pt}

\noindent
  Choose $n$ such that $2^n > |H|$.  Let $B_{nj+1}$ be the full binary
  tree of height $nj+1$.  By Proposition~\ref{prop:coloring}, if we
  label the leaves with elements of $H$, there must be at least two
  leaves at a distance $qj+1$ having the same label, for some 
  $q\leq n$.  The full binary tree $B_{qj+1}$ containing these two
  leaves evaluates to $1$, and so does  $B_{nj+1}$.  That means
  $B_{nj+1}(H)=1$ and by Lemma~\ref{lemma:passing to G^p}
  $H^{(nj+1)}=1$, so $H$ is solvable, and so is $G$.
\endproof

\bibsection


\begin{biblist}

\bib{cfp}{article}{
   author={Cannon, J. W.},
   author={Floyd, W. J.},
   author={Parry, W. R.},
   title={Introductory notes on Richard Thompson's groups},
   journal={Enseign. Math. (2)},
   volume={42},
   date={1996},
   number={3-4},
   pages={215--256},
   issn={0013-8584},
   review={\MR{1426438 (98g:20058)}},
}

\bib{geog-guzm}{article}{
   author={Geoghegan, Ross},
   author={Guzm{\'a}n, Fernando},
   title={Associativity and Thompson's group},
   conference={
      title={Topological and asymptotic aspects of group theory},
   },
   book={
      series={Contemp. Math.},
      volume={394},
      publisher={Amer. Math. Soc.},
      place={Providence, RI},
   },
   date={2006},
   pages={113--135},
   review={\MR{2216710 (2007a:20067)}},
}

\bib{levi}{article}{
   author={Levi, F. W.},
   title={Groups in which the commutator operation satisfies certain
   algebraic conditions},
   journal={J. Indian Math. Soc. (N.S.)},
   volume={6},
   date={1942},
   pages={87--97},
   review={\MR{0007417 (4,133i)}},
}


\end{biblist}
\end{document}